\documentclass[final,12pt]{amsart} 
\usepackage{amsmath, bm} 
\usepackage{amssymb} 
\usepackage{graphicx} 
\def\e{\varepsilon}

\title[INVARIANT TORI OF A TWO DIMENSIONAL PERIODIC SYSTEM] {INVARIANT TORI OF A TWO DIMENSIONAL PERIODIC SYSTEM WITH THE LINEAR-CUBIC UNPERTURBED PART
  }

\author[ V. V. Basov, V. G. Romanovski  and A. S. Zhukov ]
{  Vladimir V. Basov$^{1*}$, Valery  G. Romanovski$^{2,3,4}$, Artem S. Zhukov$^{1}$}

\address{$^1$Faculty of Mathematics and Mechanics, St. Petersburg State University, 
 Universitetsky prospekt, 28, 198504, Peterhof, St. Petersburg, Russia}

\email{vlvlbasov@rambler.ru,  artzhukov1111@gmail.ru}

\address{$^2$ Faculty of Electrical Engineering and Computer Science, University of  Maribor, SI-2000 Maribor, Slovenia}

\address{$^3$ Center for Applied Mathematics and Theoretical Physics, SI-2000 Maribor, Slovenia}
\address{$^4$Faculty of Natural Science and Mathematics, University of Maribor,
  Koro\v ska cesta 160, SI-2000 Maribor, Slovenia}

\email{valerij.romanovskij@um.si}



\begin{document}

\begin{abstract} Two classes of time-periodic systems of ordinary differential equations with a small parameter $\e\ge0,$ those with "fast"\ and "slow"\ time, are studied.
Right-hand sides of these systems are three times continuously differentiable with respect to phase variables and the parameter, the corresponding unperturbed systems are autonomous,
conservative and have nine equilibrium points.
 For the perturbed systems, which do not  depend on the parameter explicitly, we obtain the conditions yielding that the initial system has a certain number of two-dimensional invariant surfaces homeomorphic to a torus
for each sufficiently small values of parameter $\e$ and the  formulas of such surfaces. 
A class of systems with seven invariant surfaces enclosing different configurations  of equilibrium points is studied  as an example of   applications of our method.
\end{abstract}

\maketitle

\section{Introduction.}


We study a periodic two-dimensional  system of ODE's  with a small parameter $\e$ defined by the equations  
\begin{equation}\label{sv}
\dot x=\big(\gamma(y^3-y)+X(t,x,y,\e)\e\big)\e^\nu,\ \ \dot y=\big(\!-(x^3-x)+Y(t,x,y,\e)\e\big)\e^\nu,  
\end{equation}
where $\gamma\in(0,1],\ \nu=0,1;$
$X,Y$ are continuous $C_{x,\,y,\,\e}^3$ functions,  $T$-periodic in $t$ for $t\in \mathbb{R}^1,$
$|x|<M_x,\ |y|<M_y$ $(M_x>\sqrt{2},\,M_y>\sqrt{1+\gamma^{-1/2}});$ 
$\e\in[0,\e_0).$

In essence, formula \eqref{sv} determines two different systems: one with $\nu=0,$ another with $\nu=1.$ 
Comparing these systems we can say that the system with $\nu=1,$ which  is usually called the standard system, has the "fast"\ time, 
because reducing it to the system with $\nu=0$ we obtain the period $T\e.$

It is natural to refer to the autonomous system
\begin{equation}\label{snv}
\dot x=\gamma(y^3-y)\e^\nu,\quad \dot y=-(x^3-x)\e^\nu \qquad (\gamma\in(0,1],\ \nu=0,1)
\end{equation}
as to the system of the first approximation or the unperturbed system with respect to \eqref{sv}.
System \eqref{snv} is conservative. It has nine equilibrium points and its phase plane is filled, in addition to the equilibrium points, with the closed orbits
and separatrices determined by the integral \,$(x^2-1)^2+\gamma(y^2-1)^2=a.$

The goal of this paper is to find, for system \eqref{sv} with any sufficiently small $\e>0,$ a certain number of 2-dimensional cylindrical invariant surfaces homeomorphic to the torus, which is obtained by factoring time with respect to the period. The projections of such surfaces on the phase plane are contained in a small neighborhood of corresponding closed orbits of the unperturbed system \eqref{snv}.

We explicitly write out conditions  (depending  on the parameter $\gamma$) on the unperturbed functions $X(t,x,y,0)$
and  $Y(t,x,y,0)$ under which the perturbed system \eqref{sv}
has invariant surfaces described  above and obtain asymptotic expansions in powers of $\e$ 
for  each of them. 
Besides, we provide conditions that specify a class of systems with seven invariant surfaces. Phase plane projections of such surfaces encloses one, three or nine equilibrium points of the unperturbed system.

We also provide  examples  of systems \eqref{sv} for which the three obtained bifurcational equations have admissible solutions $c_i\ \ (i=0,1,2).$

Systems \eqref{sv} with $\nu=0$ and $\nu=1$ can be studied simultaneously, because the invariant tori can be found using the method 
developed in [1,\,2] and essentially modified in [3--7].  But significant differences arise in the process of averaging the systems obtained 
as the result of a special polar coordinates change. 

One of the most significant differences of this paper from the previous works is the fact that equilibrium points are located not only on the abscissa axis. Due to this, a new interesting phenomenon that has been  never observed before arises -- it is  the inability to perform the  special polar coordinate change in a small neighborhood of singular points $(1,0)$ and $(-1,0)$ despite  moving the origin to these points.

\section{Parametrization of the orbits of the unperturbed system.\\} 

\subsection{Construction of the Phase Portrait.}  
Consider an autonomous conservative system with nine equilibrium points similar to system \eqref{snv}
defined by the equations 
\begin{equation}\label{cs}
C'(\varphi)=\gamma(S^3(\varphi)-S(\varphi)),\ \ S'(\varphi)=-(C^3(\varphi)-C(\varphi))\quad (0<\gamma\le 1).
\end{equation}

For each $a>0,$ we consider the set $\Gamma_a$ of closed orbits on the plane $(C,S)$ determined by the integrals of  system \eqref{cs} given by the equation 
\begin{equation}\label{oi} 
  (C^2-1)^2+\gamma(S^2-1)^2=a. 
\end{equation}

Obviously, systems \eqref{cs} and \eqref{snv} have the same closed orbits, but if $\nu=1$ the
 solutions of \eqref{cs} and \eqref{snv} determine different trajectories.

For system \eqref{cs} five of nine singular points are centers. Points $(\pm 1,1),$ $(\pm 1,-1)$ are solutions of \eqref{oi} when $a=0,$ and
$(0,0)$ is the  solution of \eqref{oi} when $a=1+\gamma.$ The other 
 singular points are the saddle points: points $(\pm 1,0)$ are the solutions of \eqref{oi} when $a=\gamma,$
and $(0,\pm 1)$ are the solutions of \eqref{oi} when $a=1.$


It is sufficient  to describe orbits or parts of orbits that lie in the first quadrant, because equation \eqref{oi} is invariant with  respect to the change  $C \to -C$,  $S \to -S$. 
We denote such set of orbits by  $\Gamma_a^*.$ 

The extremal values of the curves that belong to $\Gamma_a^*$ are:
\begin{equation}\label{exc}
\begin{matrix} 
r_\gamma=\sqrt{1+\gamma^{1/2}},\ \ l_\gamma=\sqrt{1-\gamma^{1/2}},\hfill\\
r_i=\sqrt{1-(1-\gamma)^{1/2}},\ \ r_e=\sqrt{1+(1-\gamma)^{1/2}},\ \ u_e=\sqrt{1+\gamma^{-1/2}};\hfill\\
r_{0i}^0,l_1^0=\sqrt{1-(a-\gamma)^{1/2}},\ \ r_{0e}^0,r_1^0=\sqrt{1+(a-\gamma)^{1/2}},\ \ r_{0e}^1,r_2^1=\sqrt{1+a^{1/2}},\hfill\\ 
u_{0i}^0=\sqrt{1-((a-1)/\gamma)^{1/2}},\ u_{0e}^0=\sqrt{1+((a-1)/\gamma)^{1/2}},\ u_{0e}^1,u_1^1,u_2^1=\sqrt{1+(a/\gamma)^{1/2}},\hfill\\ 
r_1^1=\sqrt{1+a^{1/2}},\ \ l_1^1,l_2^1=\sqrt{1-a^{1/2}},\ \ lo_2^1=\sqrt{1-(a/\gamma)^{1/2}},\hfill \end{matrix}
\end{equation}
where  $r$, $l$, $u$, $lo$, $i$, $e$ means right,  left, upper, lower, internal and  external, respectively,
the first five constants characterize $\Gamma_\gamma^*$ and $\Gamma_1^*,$ 
for the other constants the superscript determines the value of the second coordinate, 
while the  subscript determines the  number of a class which we will introduce below.
\begin{figure}
\includegraphics[scale=0.28]{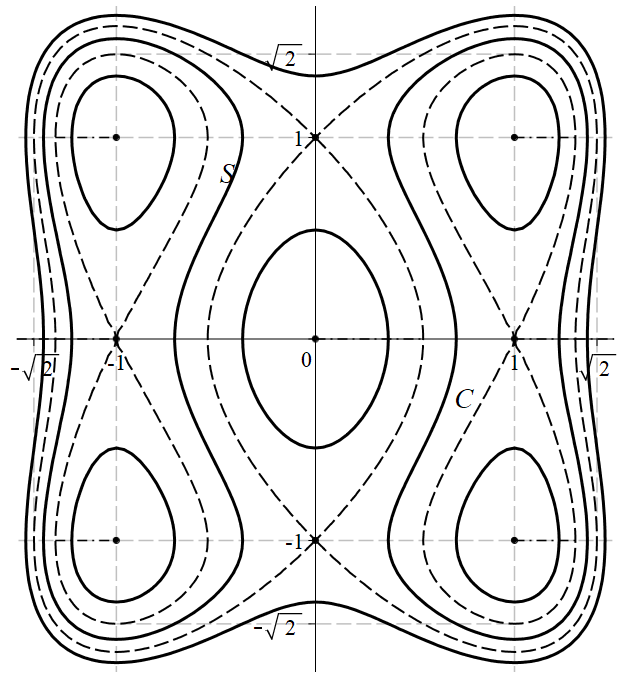} 
\includegraphics[scale=0.28]{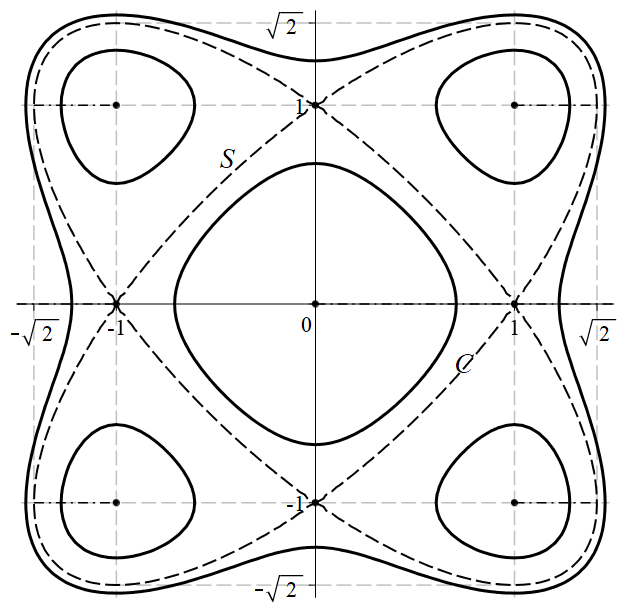}  
\caption{Phase portrait of the unperturbed system with $\gamma=1/2$ (left) and $\gamma=1$ (right).}
\end{figure}
We first  consider the "separatrix" curves determined by equation \eqref{oi}, namely, the curves  
passing  through the singular points of system \eqref{cs}.
%
%

For $\gamma<1$  the set $\Gamma_\gamma$ $(a=\gamma)$  consists of four closed curves (or two "eights"). 
Two  of those curves contact each other at the singular point $(1,0)$ and enclose points $(1,\pm 1),$ 
other two contact each others at the singular point $(-1,0)$ and enclose singular points $(-1,\pm 1).$
Thus, $\Gamma_\gamma^*$ is a top part of right "eight" and has the following extremal points: $(r_\gamma,1),\, (1,\sqrt 2),\, (l_\gamma,1),\, (1,0).$

For $\gamma<1$ the  set $\Gamma_1$ $(a=1)$ consists of two closed curves: the internal curve $\Gamma_{1i}$ and the external one $\Gamma_{1e}.$ 
These curves contact each other at the singular points $(0,1)$ and $(0,-1).$  The set $\Gamma_{1i}$ encloses the singular point $(0,0)$,  
and one of each "eights"  mentioned above is located in the two areas between $\Gamma_{1i}$ and $\Gamma_{1e}.$ 
$\Gamma_{1i}^*$ has the extremal points $(r_i,0),\,(0,1),$ and $\Gamma_{1e}^*$ has the extremal points $(r_e,0),\,(\sqrt 2,1),\,(1,u_e),\,(0,1).$ 

For $a=\gamma=1$ $r_i,r_e=1$ and $\Gamma_\gamma$ coincides with $\Gamma_1,$ 
namely, $\Gamma_{1i}^*$ and $\Gamma_{1e}^*$ contact each other at the point $(1,0)$ creating the  top part of the right "eight".

{\bf Remark 1.} \ {\sl We choose the bounds  on constants $M_x,M_y,$ which determine the domain of system \eqref{sv},  in such  way that 
all separatrices of unperturbed system \eqref{snv} lie inside of it. 
Moreover, it follows from equation \eqref{oi} and inequalities  $|x|<M_x, |y|<M_y$ that $a<M=\min\{\gamma(M_y^2-1)^2),(M_x^2-1)^2\}>1.$ }

For $a\ne \gamma1$ the  set $\Gamma_a$ consists of closed orbits of system \eqref{cs}.
The separatrices separate the   set $\Gamma_a$ into three classes, which we denote by  0],\,1],\,2]
and define as follows:

\smallskip
0]\, $a>1.$ \ For  $1<a<1+\gamma$ the set $\Gamma_a$ consists of two closed orbits: the 
inner orbit $\Gamma_{ai}$ which encloses $(0,0)$ and lies inside $\Gamma_{1i},$
and the outer  orbit $\Gamma_{ae}$ which encloses $\Gamma_{1e}.$
Then the points $(r_{0i}^0,0),$ $(0,u_{0,i}^0)$ are extremal for $\Gamma_{ai}^*,$ 
the points $(r_{0e}^0,0),$ $(r_{0e}^1,1),$ $(1,u_{0e}^1),$ $(0,u_{0e}^0)$ are extremal for $\Gamma_{ae}^*,$ 
moreover, $r_{0i}^0\in (0,r_i),$ $r_{0e}^0\in (r_e,\infty).$  Since  $\Gamma_{ai}$ degenerates into
the  point $(0,0)$ for $a=1+\gamma,$
for $a\ge 1+\gamma,$ only the closed orbit $\Gamma_{ae}$ exists.

As the  result, the class 0] naturally splits into two subclasses: 
\,$0_i]$ for closed orbits $\Gamma_{ai}$ with $1<a<1+\gamma$\, and \,$0_e]$ for closed orbits $\Gamma_{ae}$ with $a>1.$

\smallskip
1]\, $\gamma<a<1.$ The set $\Gamma_a$ consists of two closed orbits. The top half of  the right orbit $\Gamma_a^*,$ 
is located between $\Gamma_\gamma^*$ and $\Gamma_{1i}^*\cup\Gamma_{1e}^*$\, and
has the  following extremal points: $(r_1^0,0),\,(r_1^1,1),\,(1,u_1^1) $ $(l_1^1,1),\,(l_1^0,0),$ where $r_1^0\in (1,r_e).$

\smallskip
2]\, $0<a<\gamma.$ The  set $\Gamma_a$ consists of four closed orbits.
The orbit $\Gamma_a^*$ encloses $(1,1)$ and lies inside of $\Gamma_\gamma^*,$
its extremal points are $(r_2^1,1),\,(1,u_2^1),\,(l_2^1,1),\,(1,lo_2^1),$ where $r_2^1\in (1,r_\gamma).$

\subsection{Parametrization of the closed orbits.}  The
parameter $a$ does not  determine a specific closed orbit from $\Gamma_a,$ therefore it cannot be 
used  to parametrize the orbits. 
To define the parametrization  we use the  extremal points: $(r_{0i}^0,0),(r_{0e}^0,0)$ in class 0], 
\,$(\pm r_1^0,0)$ in class 1], \,$(\pm r_2^1,1)$ and  $(\pm r_2^1,-1)$ in class 2]. 
Each such point defines  the corresponding closed orbit and  the parameter $a$ from \eqref{oi} 
can be explicitly expressed through $r$ mentioned above

Thus, an arbitrary closed orbit of system \eqref{cs} is parameterized by the functions $C(\varphi),S(\varphi).$ 
These functions are solutions of the  initial value problem 
\begin{equation}\label{nd} 
\begin{matrix} C(0)= b_{kl},\ \ S(0)=l \quad (k,l=0,\pm 1,\ (k,l)\ne (0,\pm 1));\\ 
  b_{00}=\left[\begin{matrix} b_{00}^i=r_{0i}^0\in (0,r_i)\hfill\\ 
                              b_{00}^e=r_{0e}^0\in (r_e,r_M)\end{matrix}\right.,\ \ 
  \begin{matrix} b_{10}=r_1^0\in (1,r_e),\\ 
                b_{-10}=-r_1^0,\hfill\end{matrix}\ \ 															
  \begin{matrix} b_{1,\pm 1}=r_2^1\in (1,r_\gamma),\\ 
                b_{-1,\pm 1}=-r_2^1;\hfill\end{matrix} \\
  a_{b_{kl}}=(1-|l|)\gamma +(b_{kl}^2-1)^2;\ \ r_M = \sqrt{1+(M-\gamma)^{1/2}},  \end{matrix} 
\end{equation}
where  $k$ determines the "shift"\ along the  abscissa axis, $l$ determines the "shift"\ along the  ordinate axis
and $|k|+|l|$ determines the  class number related to the  parameterized closed orbit, and the   
constants $r$ are introduced in \eqref{exc}. 
%

Let us note that the restriction  $b_{00}^e<r_M$ is introduced according to Remark 1 so we can work with  system \eqref{sv}.

We will denote the real analytic $\omega_{b_{kl}}$-periodic solution of the  initial value problem with the initial conditions \eqref{nd} 
by  $(C(\varphi,b_{kl}),S(\varphi,b_{kl})).$

\subsection{ Calculation of periods.} 
Let us introduce the auxiliary functions  
\begin{equation}\label{Spm} S^{\pm}(C^2(\varphi))=\sqrt{1\pm\gamma^{-1/2}(a_{b_{kl}}-(C^2(\varphi)-1)^2)^{1/2}}.\end{equation}
Then, in the first  integral \eqref{oi}  $S=\pm S^-(C^2)$ or $S=\pm S^+(C^2),$ 
and the 
signs are chosen according to the  location of points from the closed orbit 
parameterized by the solution $(C(\varphi),S(\varphi))$ with respect to the  lines $S=0,\pm 1.$

We can write down the first equation of system \eqref{cs} as  
\begin{equation}\label{dp} d\varphi=(\gamma(S^3(\varphi)-S(\varphi)))^{-1}dC(\varphi). \end{equation}

{\bf Proposition 1.} The following equalities hold:
\begin{equation}\label{ome}
\begin{matrix} \qquad\ \, \omega_{b_{00}^i}=4\varphi_i^-,\ \ \omega_{b_{00}^e}=4(\varphi_e^+ + \varphi_e^-),\\  
  \omega_{b_{\pm 10}}=2(\varphi_l^- + \varphi_u^+ + \varphi_r^-),\ \ 
	\omega_{b_{1,\pm 1}}=\omega_{b_{-1,\pm 1}}=\varphi_2 ^- + \varphi_2^+,\end{matrix} 
\end{equation} 
where 
$\displaystyle \varphi_i^-=\int_{b_{00}^i}^0 \zeta^-dC;$ \
$\displaystyle \varphi_e^+=\int_0^{r_{0e}^1} \zeta^+dC,$ \
$\displaystyle \varphi_e^-=\int_{r_{0e}^1}^{b_{00}^e} \zeta^-dC;$ \
$\displaystyle \varphi_l^-=\int_{l_1^0}^{l_1^1} \zeta^-dC,$ \
$\displaystyle \varphi_u^+=\int_{l_1^1}^{r_1^1} \zeta^+dC,$ \ 
$\displaystyle \varphi_r^-=\int_{r_1^1}^{b_{10}} \zeta^-dC;$ \
$\displaystyle \varphi_2^-=\int_{r_2^1}^{l_2^1} \zeta^-dC,$ \
$\displaystyle \varphi_2^+=\int_{l_2^1}^{b_{11}} \zeta^+dC,$\, 
the limits of integration are defined in \eqref{exc}, with  $a=a_{b_{kl}}$ from \eqref{nd}, 
and $\zeta^{\pm}(C^2)=(\gamma({S^{\pm}}^3(C^2)-S^{\pm}(C^2)))^{-1}.$

{\it Proof}.  Let us calculate the  period of solution $(C(\varphi,b_{10}),S(\varphi,b_{10}))$
(this solution parameterizes the closed orbit from the class 1]  when $b_{10}=r_1^0\in (1,r_e)$ and $C(0)=r_1^0,\,S(0)=0$).

There exists a $\varphi_*$ such that 
 $C(\varphi_*)=l_1^0,$ $S(\varphi_*)=0$ \ $(0<\varphi_*<\omega_{b_{10}}).$
The constants $l_1^0,l_1^1,r_1^1$ from \eqref{exc} are abscissas of the
 left intersection point of the  closed orbit and the  
lines $S=0$ and $S=1,$ and $S'(0)=r_1^0-(r_1^0)^3<0.$
Therefore,  integrating  \eqref{dp} from $\varphi_*$ to $\omega_{b_{10}},$ 
we get 
$$\omega_{b_{10}} - \varphi_* = \varphi_l^- + \varphi_u^+ + \varphi_r^-.$$ 
This calculated value is the  "period of time,"\ when the  orbit is located in the first quadrant: below and above the line \,$S=1.$ 
Similarly,  integrating formula \eqref{dp} from $0$ to $\varphi_*$ and  taking into account  
that $l^1_1,r^1_1$ are also the  abscissas of the intersection points of the closed orbits and the line $S=-1,$ 
we get $\varphi_* = \varphi_r^- + \varphi_u^+ + \varphi_l^-.$ 
Hence, $\varphi_*=\omega_{b_{10}}/2$ 
and $\omega_{b_{10}}=2(\varphi_l^- + \varphi_u^+ + \varphi_r^-);$
$\omega_{b_{-10}}=\omega_{b_{10}}$ due to the symmetry.

Computations are similar for classes $0_i$], $0_e$] and 2]. Moreover, $C(\omega_{b_{10}}/4)=0,$ $S(\omega_{b_{10}})/4=u_0^0$ for  the class 0]. \ $\Box$

\section{Dynamics in a  neighborhood of a closed orbit.} 

\subsection{ Monotonicity indicator of the angular variable.} 
For each $b_{kl}$ from \eqref{nd} we  consider the function defined  on the  periodic solution $(C(\varphi,b_{kl}),S(\varphi,b_{kl}))$ of system \eqref{cs} by the formula 
\begin{equation}\label{p} 
\alpha_{kl}(\varphi)=C'(\varphi)(S(\varphi)-l)-(C(\varphi)-k)S'(\varphi). 
\end{equation}
Using formula \eqref{oi} the  function can be written in a simpler form as 
$$\alpha_{kl}=a_{b_{kl}}-1-\gamma+C^2+\gamma S^2-k(C^3-C)-\gamma l(S^3-S).$$
%

Let us study how the function $\alpha_{kl}(\varphi)$ changes its sign along the orbits which 
go counter-clockwise when  $\varphi$ increases for the class $0_i]$, and go clockwise when $\varphi$ increases for the other classes.

The passage to a neighborhood of an arbitrary closed orbit is possible only if $\alpha_{kl}(\varphi)$ is a function of a fixed sign, because geometrically the sign of  $\alpha_{kl}$
reflects the  monotonicity of the  change of  the angular variable assuming we observe the movement along the closed orbit from the point $(k,l).$

This is the reason why the same function $\alpha_{00}$ cannot be used for passing to orbits from classes 1] and 2].
Assuming we observe the movement from the origin  the  polar angle is not  changing  monotonically.
The subtraction of constants $k$ and $l$ in formula \eqref{p} means the  origin shifts 
 to the point $(k,l)$ in system \eqref{cs}.

\subsection{Monotonicity  of the angular variable in class 0].} 
Let us show that $\alpha_{00}=a_{b_{00}}-1-\gamma+C^2(\varphi)+\gamma S^2(\varphi)$ is a function of a  fixed sign. 

We have: \,$\alpha_{00}'=2C(\varphi)C(\varphi)'+2\gamma S(\varphi)S'(\varphi)=2\gamma C(\varphi)S(\varphi)(S^2(\varphi)-C^2(\varphi)).$
%

For the class $0^m]\ (a=a_{b_{00}^i}\in (1,1+\gamma)),$ orbits are located in the first quadrant when $\varphi\in[0,\omega/4]$ and pass through the points $(r_{0i}^0,0)$
when $\varphi=0$ and $(0,u_{0i}^0)$ when $\varphi=\omega/4.$ Constants $r_{0i}^0=b_{00}^i$ and $u_{0i}^0$ are given in \eqref{exc}.

Note that  $\alpha_{00}'(\varphi)<0$ when $C_0(\varphi)>S_0(\varphi)$ and $\alpha_{00}'(\varphi)>0$ when $C_0(\varphi)<S_0(\varphi).$ 
Therefore, the function $\alpha_{00}(\varphi)$ takes maximum values at the endpoints $[0,\omega/4].$
However, $\alpha_{00}(0)=a-1-\gamma+(r_{0i}^0)^2=(a-\gamma)^{1/2}((a-\gamma)^{1/2}-1)<0,$
$\alpha_{00}(\omega/4)=a-1-\gamma+\gamma (u_{0i}^0)^2=(a-1)^{1/2}((a-1)^{1/2}-\sqrt\gamma)<0.$

Thus, due to the symmetry, we conclude that  $\alpha_{00}(\varphi)<0$ for any $\varphi.$

For the class $0^p]\ (a=a_{b_{00}^e}>1),$  the orbits are located in the first quadrant when $\varphi\in[3\omega/4,\omega]$ and pass through the point $(0,u_{0e}^0)$
when $\varphi=3\omega/4$ and the point  $(r_{0e}^0,0)$ when $\varphi=\omega.$ The constants $r_{0e}^0=b_{00}^e$ and $u_{0e}^0$ are given in \eqref{exc}.

We observe also that  $\alpha_{00}'(\varphi)>0$ when $C_0(\varphi)<S_0(\varphi)$ and $\alpha_{00}'(\varphi)<0$ when $C_0(\varphi)>S_0(\varphi).$ 
Therefore, the function $\alpha_{00}(\varphi)$ takes minimum values at the endpoints $[3\omega/4,\omega].$

However, $\alpha_{00}(3\omega/4)=a-1-\gamma+\gamma (u_{0e}^0)^2>0$ and $\alpha_{00}(0)=a-1-\gamma+(r_{0e}^0)^2>0.$ 
Thus, for each $\varphi,$ we conclude that  $\alpha_{00}(\varphi)>0$ due to the symmetry.
%

\subsection{Monotonicity of the angular variable for the  class 1].} 
We first  prove that for $b_{k0}$ from $(\ref{nd})$ close to $k$ the  function $\alpha_{k0}=\alpha_{k0}(\varphi,b_{k0})$ $(k=\pm 1)$ from $(\ref{p})$ 
alternates its sign.

Let $k=1.$ Then, for $b_{10}\rightarrow 1_{+0},$ the orbits of system $(\ref{cs})$ converge to the right "eight" from the outside.
%
%
Let us show that, for each $\gamma\in(0,1),$ $\alpha_{10}(\varphi,b_{k0})$ takes both positive and negative values assuming we move along the right "eight". 

For $b_{10}=1,$ in formula \eqref{exc} $a_{b_{10}}=\gamma,$ therefore,
function $\alpha_{10}(\varphi,1)=\gamma S^2(\varphi,1)-C^3(\varphi,1)+C^2(\varphi,1)-C(\varphi,1)-1.$
It is possible to find such $\varphi_*$, that $C(\varphi_*,1)=1,$ $S(\varphi_*,1)=\sqrt{2},$ 
because $(1,\sqrt{2})$ is the  apex point of $\Gamma^*_{\gamma}.$ Then $\alpha_{10}(\varphi_*,1)=2\gamma>0.$ 

For each $\check \varphi,$ such that  $C(\check \varphi,1)\in((1-\gamma)^{1/2},1),$ $|S(\check \varphi,1)|\in(0,1),$ 
in integral \eqref{oi} $\gamma(S^2-1)^2=\gamma-(C^2-1)^2\,(\ge 0)$ or
$S^2(\check \varphi,1)=1-(1-\gamma^{-1}(C^2(\check \varphi,1)-1)^2)^{1/2}.$ 
Omitting  the arguments  $\check \varphi,1$ we obtain $\alpha_{10}=\gamma-\gamma^{1/2}(\gamma-(C^2-1)^2)^{1/2}-(C-1)^2(C+1).$ 
Then $\alpha_{10}<0 \Leftrightarrow (\gamma - (C-1)^2(C+1))^2<\gamma(\gamma-(C_1^2-1)^2) \Leftrightarrow C(C^2+\gamma-1)<0,$ 
which is the correct statement.

Now, continuity of $\alpha_{10}$ as a function of $b_{10}$ implies that for $b_{10}$ close enough to zero, the 
function $\alpha_{10}(\varphi,b_{10})$ also alternates its sign.
It changes its sign on any closed orbit from some neighborhood of the right "eight,"\ 
to be more specific, on the  part of the orbit located to the left of the  singular point $(1,0)$ of system \eqref{cs}.

The statement  that $\alpha_{-1,0}(\varphi,b_{-1,0})$ alternates its sign can be proved  similarly.

Thus, the class 1] is different from the  other  classes in the sense  that the 
 algorithm used in the paper requiring the  function $\alpha_{k0}$ to be of a fixed sign 
can be used only for  $\gamma$ which  are separated from 1 and on a  smaller interval of $b_{k0}.$ 

Thereby, we introduce  the following restriction  for the  class 1]: 
\begin{equation} \label{gamma*}
  0<\gamma\le \gamma_*=0.806,\quad kb_{k0}\in (b^-,b^+)\ \ (k=\pm 1),
\end{equation}
where \,$b^+=\{ \sqrt{1+(1-\gamma)^{1/2}}\,\hbox{ for }\gamma\in(0,0.5],\sqrt{1+(3.25-\gamma-3\gamma^2)^{1/2}/2}\hbox{ for }\gamma\in[0.5,\gamma_*]\},$ 
$b^-=\{1.01\,\hbox{ for } \gamma\in(0,0.2],\ 0.2\gamma+0.97\hbox{ for }\gamma\in[0.2,0.75],\ 0.6\gamma+0.67\hbox{ for }$ $\gamma\in[0.75,\gamma_*]\}.$ 

In particular, $(b^-,b^+)\subset (1,r_e)$ and $1.154<b_-(\gamma_*)<b^+(\gamma_*)<1.162.$

\begin{figure}
\includegraphics[scale=0.3]{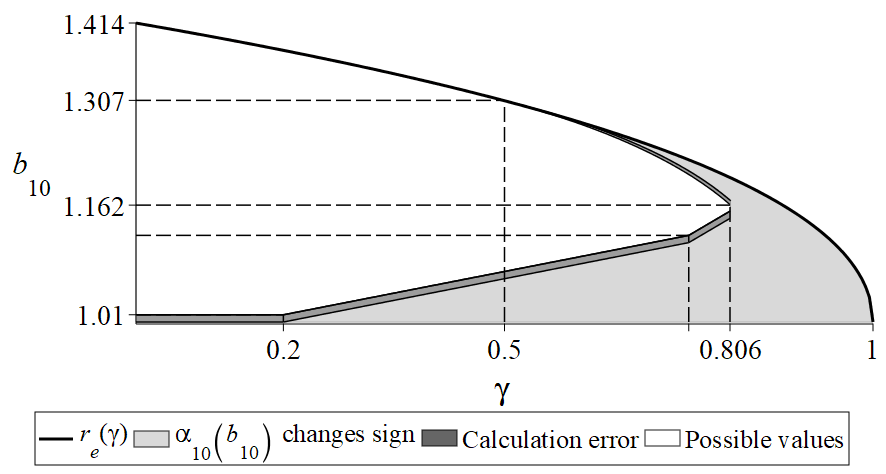}
\caption{Class 1] (k=1,l=0). Possible values of parameter $b_{10}$.}
\end{figure}

{\bf Lemma 1.}  {\it 
Let  
\,$b_*^-=\{1\,\hbox{ for } \gamma\in(0,0.2],\ 0.2\gamma+0.96\hbox{ for }\gamma\in[0.2,0.75],\ 
0.6\gamma+0.66$ $\hbox{for }\gamma\in[0.75,\gamma_*+10^{-3}] \},$
$b_*^+=\{r_e\, \hbox{ for }\gamma\in(0,0.5],$ $\sqrt{1+(3.5-\gamma-\gamma^2)^{1/2}/2}+0.02(\gamma-0.5)$ 
$\hbox{for }\gamma\in[0.5,\gamma_*+10^{-3}]\}.$

Then,\\
$1)$ for all $\gamma\in(0,\gamma_*),\, $ and  for all $ b_{k0}\!:\,kb_{k0}\in (b^-(\gamma),b^+(\gamma))\Rightarrow 
\alpha_{k0}(\varphi,b_{k0})>0;$\\
$2)$ for all  $\gamma\in (\gamma_*+10^{-3}),$ and for all  $b_{k0}\!:\,kb_{k0}\in(0,r_e)\Rightarrow 
\alpha_{k0}(\varphi,b_{k0})$ alternates its sign; \hfill\\  
$3)$ for all  $\gamma\in(0,\gamma_*+10^{-3}),\ $ 
for all  $ b_{k0}\!:\,kb_{k0}\in(1,b_*^-(\gamma))\cup(b_*^+(\gamma),r_e)$ 
$\Rightarrow \alpha_{k0}(\varphi,b_{k0})$ alternates its sign \ $(k=\pm 1).$  }

{\it Proof.} Functions and constants used in Lemma 1 were calculated approximately using  MAPLE by tracking 
the   sign of  $\alpha_{k0}$ while changing parameters $\gamma$ and $b_{k0}(\gamma)$ 
 with the step $10^{-3}.$

\subsection{Monotonicity of the angular variable in the class 2].} 
Let us now show that $\alpha_{kl}(\varphi,b_{kl})$ $(k,l=\pm1)$ from $(\ref{p})$ is positive.
%

We have that  \,$\alpha_{kl}'(\varphi,b_{kl})=\gamma  (C-k) (S-l)( kC- lS)(3kl C S+k C+ l S +1).$

Let $lk=1.$ Since  $lS>0>-( C+k)/(3 C+k),$ by multiplying the  inequality  above by $k(3C+k) ( >0)$ 
we get that $3 k l C S+k C+l S+1>0$. 

Considering that the  movement along the orbit is   clockwise the 
 analysis of  the sign of  $\alpha_{kl}$'s shows
that the function $\alpha_{kl}(\varphi,b_{kl})$ has  local minimums   at the points $\varphi_1,\varphi_2,\varphi_3$ determined by the  equations
$ C(\varphi_1)=k,$ $S(\varphi_1)=l\cdot  lo_2^1;$\, $ S(\varphi_2)=l,$ $C(\varphi_2)=k l_2^1;$\, $kC(\varphi_3)= lS(\varphi_3)>1.$ 

According to formulas $(\ref{exc})$ 
$l_2^1=\sqrt{1-a^{1/2}},\ \ lo_2^1=\sqrt{1-(a/\gamma)^{1/2}},$ hence,
$a=\gamma(1-(lo^1_2)^2)^2=(1-(l^1_2)^2)^2,$ then,
$\alpha_{kl}(\varphi_1)=a+\gamma((lo^1_2)^2 - 1)(1-lo^1_2 )=\gamma lo^1_2(lo^1_2 - 1)^2(lo^1_2+1)>0.$
Similarly, $\!\alpha_{kl}(\varphi_2)=a+((l^1_2)^2 - 1)(1-l^1_2 )=l^1_2(l^1_2 - 1)^2(l^1_2+1)>0.$ 
Finally, $ \alpha_{kl}(\varphi_3)=(1+\gamma)(S^2 - 1)S(S-l)>0.$ 

The property that  $\alpha_{kl}>0$ when $kl=-1$ can be proved analogously.

\subsection{Special polar coordinates.}   We introduce the notations: 
\begin{equation}\label{pf} 
\begin{matrix} 
	\mu_{kl} = ((b_{kl}+k)^2 - 1)^2;\ \ C_k(\varphi)=C(\varphi)-k,\ S_l(\varphi)=S(\varphi)-l; \hfill \\
	R_{kl}(t,\varphi,r,\e)=C'Y(\bar \kappa_{kl})-S'X(\bar \kappa_{kl}),\ \ (\bar \kappa_{kl})=(t,C+C_kr,S+S_lr,\e),\hfill\\ 
  \Phi_{kl}(t,\varphi,r,\e)=\alpha_i^{-1}(1+r)^{-1}(S_lX(\bar \kappa_{kl})-C_kY(\bar \kappa_{kl}));\hfill\\ 
	R_{kl}^o=R_{kl}(t,\varphi,0,0)=C'Y(t,C,S,0)-S'X(t,C,S,0),\hfill\\ 
	\Phi_{kl}^o=\Phi_{kl}(t,\varphi,0,0)=\alpha_{kl}^{-1}(S_lX(t,C,S,0)-C_kY(t,C,S,0)).\hfill\end{matrix}   
\end{equation}
Then \,$R_{kl}=R_{kl}^o+(\!{R_{kl}}'_r)^or+(\!{R_{kl}}'_\e)^o\e+O((|r|+\e)^2),$ 
 $\Phi_{kl}=\Phi_{kl}^o+O(|r|+\e),$ 
where \,$(\!{R_{kl}}'_r)^o=C'(C_kY_x'(t,C,S,0)+S_lY_y'(t,C,S,0))-
                         S'(C_kX_x'(t,C,S,0)+S_lX_y'(t,C,S,0)),$ 
	$(\!{R_{kl}}'_\e)^o= $ \\
	$C'Y_\e'(t,C,S,0)-S'X_\e'(t,C,S,0).$

Let us perform the special affine-polar coordinate change in system \eqref{sv} setting 
\begin{equation}\label{pz} 
  x=C(\varphi)+C_k(\varphi)r,\ \ y=S(\varphi)+S_l(\varphi)r\quad (|r|<r_0\le 1,\ k,l=0,\pm1),
\end{equation}
where 
$(C(\varphi),\,S(\varphi))$ are 
$\omega_{b_{kl}}$-periodic solution of system \eqref{cs} 
with the initial values $C(0)= b_{kl},\ S(0)=l$ 
from \eqref{nd} and restriction \eqref{gamma*}. 

Let us differentiate \eqref{pz} with respect to $t.$ Using  \eqref{p}  we solve the obtained equalities with respect to \,$\dot r$ and $\dot \varphi.$\, 
As the  result we get
$$\alpha_{kl}(\varphi)\dot r=C'(\varphi)\dot y-S'(\varphi)\dot x,\quad 
  (r+1)\alpha_{kl}(\varphi)\dot \varphi=S_l(\varphi)\dot x-C_k(\varphi)\dot y.$$

Substituting  the right-hand sides of system $(\ref{sv})$ into the formulas given  above 
and using notations \eqref{pf}  we obtain  the  system
\begin{equation}\label{ps}
\begin{cases}\e^{-\nu}\alpha_{kl}\dot r=-\alpha_{kl}'r-p_{kl}r^2+O(|r|^3)+ \\
  \qquad\quad\quad\ \ +(R_{kl}^o+(\!{R_{kl}}'_r)^or+(\!{R_{kl}}'_\e)^o\e+O((|r|+\e)^2))\e,  \\
\e^{-\nu}\dot\varphi=1+\alpha_{kl}q_{kl}r+O(|r|^2)+(\Phi_{kl}^o+O(|r|+\e))\e, 
\end{cases}
\end{equation}
where 
\,$p_{kl}(\varphi)=3\gamma CS((S^2-1)(C-k)^2 - (C^2-1)(S-l)^2)$ \
and  \,$q_{kl}(\varphi)=\alpha_{kl}^{-2}( (C-k)^3 (2C+k) + \gamma (S-l)^3 (2S+l)).$ 

\section{Radial averaging and the determining equation.} 

\subsection{ Primary radial averaging.} 
Let us introduce the functions
$$\beta_{kl}(\varphi)=\int_{0}^{\varphi} \xi_{kl}(s)\,ds,\quad  
\xi_{kl}(\varphi)=\alpha_{kl}^{-1}(\varphi)(\alpha_{kl}'(\varphi) q_{kl}(\varphi)-\alpha_{kl}^{-1}(\varphi) p_{kl}(\varphi)).$$ 

{\bf Lemma 2.} \ {\it The function $\beta_{kl}(\varphi,b_{kl})$ is $\omega_{b_{kl}}$-periodic in $\varphi.$ }

{\it Proof.} 
We will use  formulas \eqref{ps},\eqref{oi},\eqref{p} and \eqref{pf}.

For the class 0] \,$(k,l=0)$ we have: 
\,$p_{00}=3\alpha_{00}'/2,$ $q_{00}=4\alpha_{00}^{-1}+2\mu_{00}\alpha_{00}^{-2}.$ 
Thus,  \,$\xi_{00}=\alpha_{00}^{-2}\alpha_{00}'(5/2+2\mu_{00}\alpha_{00}^{-1}),$ 
\,$\displaystyle \beta_{00}=
\!\int_{0}^{\varphi}\! (5\alpha_{00}^{-2}(s)/2+2\mu_{00} \alpha_{00}^{-3}(s))d\alpha_{00}(s)=
5(\alpha_{00}^{-1}(0)-\alpha_{00}^{-1}(\varphi))/2+\mu_{00}(\alpha_{00}^{-2}(0)-\alpha_{00}^{-2}(\varphi))=
\beta_{00}(\alpha_{00}(\varphi,b_{00})).$ 

For the class 1] \,$(k=\pm 1,\,l=0),$\, we have: 
\,$p_{k0}=\alpha_{k0}'-\gamma C_k S_l( C_k^2+2 S_l^2+3(2-k) C_k S_l^2),$
$q_{k0}=5\alpha_{k0}^{-1}+\alpha_{k0}^{-2}(2 C_k^2-\gamma S_l^2-3\mu_{k0}).$
Therefore  \,$\xi_{k0}=\alpha_{k0}^{-2}\alpha_{k0}'(4+\alpha_{k0}^{-1}(2C_k^2-\gamma S_l^2-3\mu_{k0}))-
\gamma \alpha_{k0}^{-2}C_kS_l(C_k^2+2S_l^2+3(2-k)C_kS_l^2).$

For the class 2] \,$(k,l=\pm 1),$\, we have: 
\,$p_{kl}=\alpha_{kl}'+\gamma C_kS_l(kC_k-lS_l)(3klC_kS_l+2kC_k+2lS_l),$ 
\,$q_{kl}=5 \alpha_{kl}^{-1}+\alpha_{kl}^{-2}(2 C_k^2+2\gamma S_l^2-3\mu_{kl}),$
then \,$\xi_{kl}=\alpha_{kl}^{-2}\alpha_{kl}'(4+\alpha_{kl}^{-1}(2C_k^2+2\gamma S_l^2-3\mu_{kl}))-
 \gamma \alpha_{kl}^{-2}C_kS_l(kC_k-lS_l)(2kC_k+3klC_kS_l+2lS_l).$

For the cases $k\neq0$ let us describe the  movement along the  closed orbit when $\varphi$ changes from zero to $\omega_{b_{kl}}.$
To this end,  we  express $S_l(\varphi)$ as a  composite function  using
formula \eqref{oi} for the integral, notations \eqref{pf}, constants from \eqref{exc} and the fact that $S'(0)>0$ for  system \eqref{cs}.
Thus,  we obtain:
$$\begin{matrix}
S_0(C_k(\varphi))=\{ 
	S_0^- \hbox{ for } kC_k\searrow_{l_1^1-k}^{l_1^0-k},\,
  S_0^+ \hbox{ for } kC_k\nearrow_{l_1^1-k}^{r_1^1-k},\,
  S_0^- \hbox{ for } kC_k\searrow_{b_{k0}-k}^{r_1^1-k},\hfill \\ 
	\hfill -S_0^- \hbox{ for } kC_k\nearrow_{b_{k0}-k}^{r_1^1-k},\,
         -S_0^+ \hbox{ for } kC_k\searrow_{l_1^1-k}^{r_1^1-k},\, 
				 -S_0^- \hbox{ for } kC_k\nearrow_{l_1^1-k}^{l_1^0-k} \};\\ 
S_1(C_k(\varphi))=\{ 
  S_1^+ \hbox{ for } kC_k\nearrow_{l_2^1-k}^{b_{k1}-k},\,  
	S_1^- \hbox{ for } kC_k\searrow_{l_2^1-k}^{b_{k1}-k}\};\hfill \\ 
S_{-1}(C_k(\varphi))=\{ 
  -S_1^- \hbox{ for } kC_k\nearrow_{l_2^1-k}^{b_{k,-1}-k},\, 
	-S_1^+ \hbox{ for } kC_k\searrow_{l_2^1-k}^{b_{k,-1}-k}\};\hfill \\
S_l^{\pm}=S_l^{\pm}(C_k(\varphi))=\sqrt{(1\pm\gamma^{-1/2}(a_{b_{kl}}-C_k^2(\varphi)(C_k(\varphi)+2)^2)^{1/2})^{1/2}}-l.
	\end{matrix} $$

Therefore, $\xi_{kl}=\xi_{kl}(C_k(\varphi)),$ 
because  $\alpha_{kl}(\varphi)=\alpha_{kl}(C(\varphi))=\alpha_{kl}(C_k(\varphi)+k)$ in \eqref{p}.
Thus, taking into account \eqref{dp}, $\beta_{kl}=\beta_{kl}(C_k(\varphi)).\ \ \Box$

We now  show that the  $\omega_{b_{kl}}$-periodic averaging change
\begin{equation}\label{pave}
r= \alpha_{kl}^{-1}(\varphi,b_{kl})(z+ \beta_{kl}(\varphi,b_{kl})z^2) \quad (k,l=0,\pm1 ),
\end{equation}
transforms  system $(\ref{ps})$ into the system 
\begin{equation}\label{fix}
\begin{cases} \e^{-\nu}\dot z=O(|z|^3)+(R_{kl}^o+Z_{kl}z+(\!{R_{kl}}'_\e)^o\e+O((|z|+\e)^2))\e\\ 
	\e^{-\nu}\dot\varphi=1+q_{kl}z+O(|z|^2)+(\Phi_{kl}^o+O(|z|+\e))\e\hfill \end{cases},
\end{equation}
where \,$Z_{kl}(t,\varphi)= \alpha_{kl}^{-1}({R_{kl}}'_r)^o-2 \beta_{kl} R_{kl}^o+ \alpha_{kl}^{-1} \alpha_{kl}' \Phi_{kl}^o.$

Differentiating  $(\ref{pave})$ with respect to systems $(\ref{ps})$ and $(\ref{fix})$ we obtain
\,$-\alpha_{kl}^{-1} \alpha_{kl}'(z+ \beta_{kl}z^2)- \alpha_{kl}^{-2} p_{kl}z^2+O(|z|^3)+
( R_{kl}^o+ \alpha_{kl}^{-1} (\!{R_{kl}}'_r)^oz+ (\!{R_{kl}}'_\e)^o\e+O((|z|+\e)^2))\e=
(1+2 \beta_{kl}z)( R_{kl}^o+ Z_{kl}z+ (\!{R_{kl}}'_\e)^o\e)\e+
( \beta_{kl}'z^2- \alpha_{kl}^{-1} \alpha_{kl}'z- \alpha_{kl}^{-1} \alpha_{kl}' \beta_{kl} z^2)(1+ q_{kl}z+ \Phi_{kl}^o\e).$

Equating the coefficients of $z\e$ and $z^2$ we get  the  formula for $Z_{kl}$ 
and the equation $\beta'_{kl}=\alpha_{kl}^{-1}(\alpha_{kl}' q_{kl}-\alpha_{kl}^{-1} p_{kl}).$ 
Solving  this equation we find   $\beta_{kl}(\varphi,b_{kl}).$ 

\subsection{ Decompositions of two-periodic functions, the Siegel condition.} \
We continue the  search for  two-dimensional invariant surfaces of system \eqref{sv}.
For this purpose  we average the  functions $R_{kl}^o,\,Z_{kl}$ and $(\!{R_{kl}}'_\e)^o$ in system \eqref{fix}.  
However for this system  $\dot\varphi=1+\ldots$ when $\nu=0,$ $\dot\varphi=\e+\ldots$ for  $\nu=1.$ 
Therefore the  following   averaging changes and their existence conditions will be different for each $\nu.$ 
To distinguish the cases the  functions, constants and formulas that depend on $\nu,$ 
will have the superscript 0 or 1.

For continuous, $T$-periodic in $t,$ real analytic and $\omega$-periodic in $\varphi$ 
functions $v=\upsilon^\nu(t,\varphi)$ we use the following decomposition depending on the value of parameter $\nu:$
$$\upsilon^\nu(t,\varphi)=\overline\upsilon^\nu+\hat\upsilon^\nu(\varphi)+
  \tilde\upsilon^\nu(t,\varphi)\qquad (\nu=0,1),$$
where $\displaystyle \overline\upsilon^\nu=\frac{1}{\omega T}
  \int_0^\omega\int_0^T\upsilon^\nu(t,\varphi)\,dt\,d\varphi$ is the average value of the  function $\upsilon^\nu,$
$\displaystyle \hat\upsilon^0=0,\ \ \hat\upsilon^1=\frac{1}{T}\int_0^T\upsilon^1(t,\varphi)\,dt-\overline\upsilon^1.$ 

Then, the  function $\tilde\upsilon^\nu(t,\varphi)$ has  zero average value with respect to $t,$ 
which implies periodicity of the function $\displaystyle\int_{t_0}^t \tilde\upsilon^\nu(\tau,\varphi)\,d\tau,$ 
which also has  zero average value by the virtue of the choice of the constant $t_0\in[0,T].$
%

\smallskip
{\bf Proposition 2.} {\it 
Assume $\nu=0$ and the periods $T$ and $\omega$ of continuous, $T$-periodic in $t,$ 
real analytic and $\omega$-periodic in $\varphi$ function $\tilde\upsilon^0(t,\varphi)$ 
satisfy the Siegel condition 
\begin{equation}\label{zig}
|pT-q\omega|>K(p+q)^{-\tau}\quad (K>0,\ \tau\ge 1,\ \ p,\,q\ - \hbox{\sl positive integers}).
\end{equation}
Then the equation
$$\dot{\tilde\zeta}^0(t,\varphi)+\tilde\zeta^{0'}(t,\varphi)=\tilde\upsilon^0(t,\varphi)$$
has the unique  solution $\tilde\zeta^0(t,\varphi),$ which has the same properties as the function $\tilde\upsilon(t,\varphi).$\footnote{See Lemma B.5 in [9,\,p.\,17].}} 

Here and in what follows we will denote the  derivative with respect to $t$ of any function which has $t$ and $\varphi$ as its arguments by  a dot, and the derivative with respect to $\varphi$ by 
a prime.


\subsection{ The choice of generating orbits.} 
 The  function 
 \begin{equation}\label{rio} 
R_{kl}^o(t,\varphi)=C'(\varphi)Y(t,C(\varphi),S(\varphi),0)-S'(\varphi)X(t,C(\varphi),S(\varphi),0),
\end{equation}
introduced in 
  \eqref{pf} will  play the key role
in our analysis. 
In \eqref{rio} \,$(C(\varphi),S(\varphi))$ is a real analytic $\omega_{b_{kl}}$-periodic solution of
the  initial value problem for   
system $(\ref{cs})$ with the initial values $C(0)=b_{kl},\,S(0)=l,$ the  parameter $b_{kl}$ is any number from \eqref{nd}
and the  period $\omega_{b_{kl}}$ is calculated in \eqref{ome}. 


Using the mentioned above decomposition, we write out $R^{o}_{kl}(t,\varphi)$ as the following sum:
$$R_{kl}^o=\overline {R_{kl}^o}+\widehat R_{kl}^o(\varphi)+\widetilde R_{kl}^o(t,\varphi),$$
where $\displaystyle \overline {R_{kl}^o}(b_{kl})={1\over T\omega}\int_0^{\omega}\int_0^T R_{kl}^o(t,\varphi)dt\,d\varphi$ 
is  the average value of $R_{kl}^o$ \,$(\omega=\omega_{b_{kl}}).$ 

For each $k,l$ $(k,l=0,\pm 1,\ (k,l)\ne (0,\pm 1)),$  we also introduce the  generating equation
\begin{equation}\label{pu} \overline {R_{kl}^o}(b_{kl})=0.\end{equation} 

Moreover, if for $\nu=0$ the periods $T$ and $\omega_{kl}$ satisfy Siegel condition \eqref{zig}, then by Proposition 2 the equation
\begin{equation}\label{gn}
\dot{\tilde g}_{kl}^0(t,\varphi)+\tilde g_{kl}^{0'}(t,\varphi)=\widetilde R_{kl}^o(t,\varphi)
\dot {\check g}_i^0+\check g_i^{0'}=R_{i0}(t,\varphi)\quad (\overline {R_{i0}}=0)
\end{equation}
 has the unique  solution $\tilde g_{kl}^0(t,\varphi).$ 
This solution has the same properties as the function $\widetilde R_{kl}^o(t,\varphi).$ 

Now we can formulate the dissipativity condition as 
\begin{equation} \label{ud}  L_{kl}^\nu\ne 0\qquad (\nu=0,1), \end{equation}
where $L_{kl}^0=\overline{Z_{kl}}-\overline {{g_{kl}^0}'q_{kl}},$\, 
 \,$L_{kl}^1=\overline{Z_{kl}}-\overline {\widehat R_{kl}^o q_{kl}};$
and $q_{kl}(\varphi)$ is from \eqref{ps}, $Z_{kl}(t,\varphi)$ is from \eqref{fix}.


{\bf Definition 1.} \ {\it The parameter $b_{kl}$ from \eqref{nd} is called admissible for \eqref{sv} 
and is denoted by $b_{kl}^*,$ if it is a solution of generating equation \eqref{pu},
condition \eqref{gamma*} holds for  $|k|+|l|=1,$
the periods \,$T$ and \,$\omega^*=\omega_{b_{kl}^*}$ satisfy   Siegel condition \eqref{zig} for $\nu=0$
and the dissipativity condition \eqref{ud} holds.}

{\bf Remark 2.} \ {\sl 
We will  be interested in systems \eqref{sv} with a nonempty set of admissible values.
For each  $b_{kl}^*,$ it will be proved that the perturbed system \eqref{sv} retains a two-periodic invariant surface homeomorphic to a torus, which is obtained by factoring the time with
respect to the period in a small  neighborhood of the cylindric surface, whose generatrix is a closed orbit of unperturbed system \eqref{snv} that passes through the point $(b_{kl}^*,l)$ 
for each sufficiently small $\e.$}

From now on  we fix an  admissible parameter $b_{kl}^*$ of system \eqref{sv} which in turn fixes
 initial values for the 
$\omega^*$-periodic initial value problem solution $(C(\varphi),S(\varphi))$ of system \eqref{cs}, which parameterizes a specific closed orbit related to the class $|k|+|l|.$ 

All subsequent changes and introduced functions are fixed by the choice of $b_{kl}^*.$ 
In particular, conditions \eqref{gamma*},\,\eqref{zig},\,\eqref{ud} hold by the  definition and 
$R_{kl}^o,$ which is a part of systems \eqref{ps} and \eqref{fix},  has  zero average value, that is,
$$R_{kl}^o(t,\varphi)=\widetilde R_{kl}^o(t,\varphi).$$

\section{The construction of invariant surfaces. The results.} 

\subsection{ Secondary radial averaging.} 
We introduce the  functions: 
$$g_{kl}^0=\overline {g_{kl}^0}+\tilde g_{kl}^0,\ \ \tilde g_{kl}^0 \hbox{ from (\ref{gn}), } \
\overline {g_{kl}^0}=(\overline {{g_{kl}^0}'\Phi_{kl}^o}-\overline {(\!{R_{kl}}'_\e)^o}-
  \overline {\tilde g_{kl}^0(Z_{kl}-{g_{kl}^0}'q_{kl})})/L_{kl}^0;$$ 
\begin{equation}\label{af}
\hat g_{kl}^{1'}=\widehat R_{kl}^o,\ \ 
\overline {g_{kl}^1}=(\overline {\widehat R_{kl}^o\Phi_{kl}^o}-\overline {(\!{R_{kl}}'_\e)^o}-
\overline {\hat g_{kl}^1(Z_{kl}-\widehat R_{kl}^oq_{kl})})/L_{kl}^1.
\end{equation}

Let us show that the coordinate change
\begin{equation}\label{save}
z=u+G_{kl}^\nu(t,\varphi,\e)\e+H_{kl}^\nu(t,\varphi,\e)u\e+F_{kl}^\nu(t,\varphi,\e)\e^2\quad (\nu=0,1),
\end{equation}
where  $G_{kl}^0=g_{kl}^0(t,\varphi),$ $H_{kl}^0=\tilde h^0_{kl}(t,\varphi),$ $F_{kl}^0=\tilde f^0_{kl}(t,\varphi),$
$G_{kl}^1=\overline {g^1_{kl}}+\hat g^1_{kl}(\varphi)+\tilde g^1_{kl}(t,\varphi)\e,$
$H_{kl}^1=\hat h^1_{kl}(\varphi)+\tilde h^1_{kl}(t,\varphi)\e,$
$F_{kl}^1=\hat f^1_{kl}(\varphi)+\tilde f^1_{kl}(t,\varphi)\e,$
transforms system \eqref{fix} into the system 
\begin{equation}\label{savesys}
\dot u=(L_{kl}^\nu u\e+O((|u|+\e)^3))\e^\nu,\ \
\dot\varphi=(1+\Theta^{\nu}_{kl}\e+q_{kl}u+O((|u|+\e)^2))\e^\nu,
\end{equation}
where, obviously, $\Theta^0_{kl}(t,\varphi)=\Phi_{kl}^o+q_{kl} g_{kl}^0,$\, 
$\Theta^1_{kl}(t,\varphi)=\Phi_{kl}^o+q_{kl}(\overline {g_{kl}^1}+\hat g_{kl}^1).$ 

We achieve this by  differentiating coordinate  change (\ref{save}) with respect to systems $(\ref{fix})$ and (\ref{savesys}) and reducing the result by $\e^\nu.$ 
The calculations yield the  identity
\begin{multline*}
(R_{kl}^o+(u+G_{kl}^\nu\e)Z_i+((\!{R_{kl}}'_\e)^o\e)\e+O((|u|+\e)^3)=\\
  L_{kl}^\nu u\e+G_{kl}^{\nu'}\e(1+\Theta_{kl}^\nu\e+q_{kl}u)+H_{kl}^{\nu'}u\e+F_{kl}^{\nu'}\e^2+
(\dot G_{kl}^\nu+\dot H_{kl}^\nu u+\dot F_{kl}^\nu \e)\e^{1-\nu}.
\end{multline*}

Let $\nu$ be equal to $0.$
Then, the coefficients of 
$\e$ constitute the equation \eqref{gn}, 
which has a single two-periodic solution $\tilde g_{kl}^0.$ 

The coefficients of $u\e$ constitute the
 equation $\dot{\tilde h}_{kl}^0+\tilde h_{kl}^{0'}=Z_{kl}-{g_{kl}^0}'q_{kl}-L_{kl}^0,$ 
which, by Proposition \,2, has a single two-periodic solution $\tilde h_{kl}^0(t,\varphi),$
since the  right-hand side of the  equation has the zero average value by the virtue of choice of  constant $L_{kl}^0$ in \eqref{ud}. 

The equation $g_{kl}^0Z_{kl}+(\!{R_{kl}}'_\e)^o={g_{kl}^0}'(\Phi_{kl}^o+g_{kl}^0q_{kl})+\tilde f_{kl}^{0'}+\dot {\tilde f}_{kl}^0$ 
is constituted by the  coefficients of $\e^2.$ 
Substituting into this equation  $Z_{kl}$ from the previous  equation we get
$\dot {\tilde f}_{kl}^0+\tilde f_{kl}^{0'}=(\overline {g_{kl}^0}+\tilde g_{kl}^0)(L_{kl}^0+{h_{kl}^0}'+\dot h_{kl}^0)+
 (\!{R_{kl}}'_\e)^o-{g_{kl}^0}'\Phi_{kl}^o.$ 
 The 
two-periodic function $\tilde f_{kl}^0(t,\varphi)$ can be found explicitly by the virtue of the  choice of constant $\overline{g_{kl}^0}$ from $(\ref{af}^0).$ 

Now, consider the case  $\nu=1.$ Let us substitute decompositions of functions $G^1,H^1$ and $F^1,$ introduced in (\ref{save}) into the  identity.

For $\e,$ we obtain the equation $R_{kl}^o=\hat g_{kl}^{1'}+\dot{\tilde g}_{kl}^1,$ 
which can be explicitly solved by splitting it into two equations:
$\hat g_{kl}^{1'}=\widehat R_{kl}^o$ and $\dot{\tilde g}_{kl}^1=\widetilde R_{kl}^o.$

For $u\e,$ we have $Z_{kl}=L_{kl}^1+\hat g_{kl}^{1'}q_{kl}+\hat h_{kl}^{1'}+\dot{\tilde h}_{kl}^1$ 
with $L_{kl}^1=\overline {Z_{kl}}-\overline{\hat g_{kl}^{1'}q_{kl}}$ from \eqref{ud}.
We express $\hat h_{kl}^1,\,\tilde h_{kl}^1$ from equations
$\hat h_{kl}^{1'}=\widehat Z_{kl}-\hat g_{kl}^{1'}q_{kl}$ and $\dot{\tilde h}_{kl}^1=\widetilde Z_{kl}$ $(\widetilde q_{kl}=0).$

For $\e^2,$ we have $(\overline {g_{kl}^1}+\hat g_{kl}^1)Z_{kl}+(\!{R_{kl}}'_\e)^o=\tilde g_{kl}^{1'}+
  \hat g_{kl}^{1'}(\Phi_{kl}^o+(\overline {g_{kl}^1}+\hat g_{kl}^1)q_{kl})+\hat f_{kl}^{1'}+\dot{\tilde f}_{kl}^1.$ 
Similarly to the case  $\nu=0,$ we obtain the equation 
$\overline {g_{kl}^1}L_{kl}^1=\overline {\hat g_{kl}^{1'}\Phi_{kl}^o}-\overline {(\!{R_{kl}}'_\e)^o}-
\overline {\hat g_{kl}^1(\hat h_{kl}^{0'} +\dot{\tilde h}_{kl}^0}).$ 
This equation allows to fix the value of $\overline {g_{kl}^1}$ introduced in $(\ref{af}^1).$

After that we write out and explicitly solve the equations for $\hat f_{kl}^{1'}$ and $\dot{\tilde f}_{kl}^1.$
$\square$


\subsection{ Angular averaging and Hale's Lemma.} 
Let us average $\Theta_{kl}^\nu(t,\varphi)$ in system \eqref{savesys}  performing the 
$T$- and $\omega^*$-periodic change of the angular variable
\begin{equation}\label{phi}
\varphi=\psi+\Delta_{kl}^\nu(t,\psi,\e)\e,
\end{equation}
where $\Delta_{kl}^0=\tilde\delta_{kl}^0(t,\psi),$ $\Delta_{kl}^1=\hat\delta_{kl}^1(\psi)+\tilde\delta_{kl}^1(t,\psi)\e.$
This change transforms (\ref{savesys}) into the system
\begin{equation}\label{phisys}
\dot u=(L_{kl}^\nu u\e+O((|u|+\e)^3))\e^\nu,\ \
\dot\psi=(1+\overline {\Theta_{kl}^\nu}\e+q_{kl}(\psi)u+O((|u|+\e)^2))\e^\nu.
\end{equation} 

It is obvious that the function $\Delta_{kl}^\nu$ from $(\ref{phi})$ is uniquely determined by the equations
$\tilde \delta_{kl}^{0'}+\dot{\tilde\delta}_{kl}^0=\Theta_{kl}^0-\overline {\Theta_{kl}^0},$ \
$\hat \delta_{kl}^{1'}=\widehat\Theta_{kl}^1,$ \ $\dot{\tilde\delta}_{kl}^1=\widetilde\Theta_{kl}^1.$

It is easy to see that the change inverse to $(\ref{phi})$ can be written in the form
$$\psi=\varphi+\Omega_{kl}^\nu(t,\varphi,\e)\e,$$
where  $\Omega_{kl}^0=-\tilde\delta_{kl}^0(t,\varphi)+\tilde\delta_{kl}^0(t,\varphi)\tilde\delta_{kl}^{0'}(t,\varphi)\e+O(\e^2),$
$\Omega_{kl}^1=-\hat\delta_{kl}^1(\varphi)+(\hat\delta_{kl}^1(\varphi)\hat\delta_{kl}^{1'}(\varphi)-
\tilde\delta_{kl}^1(t,\varphi))\e+O(\e^2)$ are $T$-periodic in $t$ and
real analytic $\omega^*$-periodic in $\varphi$.

We make the scaling change
\begin{equation}\label{scale}
u=v\e^{3/2},
\end{equation}
which transforms  $(\ref{phisys})$ into the system
\begin{equation}\label{scsys}
\dot v=( L_{kl}^\nu v\e+V_{kl}^\nu(t,\psi,v,\e)\e^{3/2})\e^\nu,\ \
\dot\psi=(1+\overline{\Theta_{kl}^\nu}\e+\Psi_{kl}(t,\psi,v,\e)\e^{3/2})\e^\nu,
\end{equation}
where $V_{kl}^\nu,\,\Psi_{kl}^\nu$  are continuous functions of their arguments in a small neighborhood of $v$ and $\e,$ 
which are continuously differentiable with respect to $v$ and $\psi,$ $T$-periodic in $t$ and $\omega^*$-periodic in $\psi.$

Indeed, $V_{kl}^\nu(t,\psi,v,\e)=O((|v|\e^{3/2}+\e)^3)\e^{-3},$
$\Psi_{kl}^\nu(t,\psi,v,\e)=q_{kl}v+O((|v|\e^{3/2}+\e)^2)\e^{-3/2},$
and the functions $O(\ldots)$ are real analytic for each $\psi$ and
three times continuously differentiable in a small neighborhood of the point $v=\e=0.$
Therefore $V_{kl}^\nu, \Psi_{kl}^\nu$ are continuously differentiable with respect to $v$ at this point.

System $(\ref{scsys})$ satisfies conditions of Hale's lemmas 2.1 and 2.2 in [8], 
hence for all sufficiently small $\e,$ it has an invariant surface of the form
\begin{equation}\label{Gamma} v=\Gamma_{kl}^\nu(t,\psi,\e)\e^{1/2}, \end{equation}
where \,$\Gamma_{kl}^\nu$ is continuous continuously differentiable, $T$-periodic in $t,$ 
and $\omega^*$-periodic in $\psi$.

\subsection{ Summary of the results.} 

In the following statements we summarize the obtained results.

{\bf Lemma 3.} \ {\it 
For each admissible $b_{kl}^*$, for each sufficiently small  \,$\e>0,$ system \eqref{fix} has the continuous continuously differentiable, 
$T$-periodic in $t,$ and $\omega^*$-periodic in $\varphi$ invariant surface
\begin{equation}\label{surz} 	z=Q_{kl}^\nu(t,\varphi,\e)\quad (\nu=0,1) \end{equation}
where $Q_{kl}^\nu=G_{kl}^\nu(t,\varphi,\e)\e+(F_{kl}^\nu(t,\varphi,\e)+
 \Gamma_{kl}^\nu(t,\varphi+\Omega_{kl}^\nu(t,\varphi,\e)\e,\e)(1+H_{kl}^\nu(t,\varphi,\e)\e))\e^2.$ }

The surface $(\ref{surz})$ is obtained by substituting the invariant surface $(\ref{Gamma})$ 
into the composition of changes $(\ref{save}),$ $(\ref{scale})$ and the change inverse to $(\ref{phi}).$

{\bf Corollary 1.}  {\it For each sufficiently small $\e>0,$ system $(\ref{ps})$ has the continuous continuously differentiable $T$-periodic in $t,$ and $\omega^*$-periodic in $\varphi$ invariant surface
\begin{equation}\label{surr} 
r=\Upsilon_{kl}^\nu(t,\varphi,\e) \end{equation}
where $\Upsilon_{kl}^\nu=\alpha_{kl}^{-1}(\varphi)(Q_{kl}^\nu(t,\varphi,\e)+
  \beta_{kl}(\varphi)(Q_{kl}^\nu(t,\varphi,\e))^2).$ }

The surface $(\ref{surr})$ is obtained by substituting the invariant surface $(\ref{surz})$ into the change $(\ref{pave}).$

{\bf Corollary 2.} \ {\it	The invariant surface $(\ref{surr})$ 
has the following asymptotic expansion: }
$$\begin{matrix} 
  \Upsilon_{kl}^0=\alpha_{kl}^{-1}(\varphi)(g_{kl}^0(t,\varphi)\e+(\tilde f_{kl}^0(t,\varphi)+\Gamma_{kl}^0(t,\varphi,0)+
	\beta_{kl}(\varphi)g_{kl}^0(t,\varphi)^2)\e^2)\!+\!O(\e^3),\\
\Upsilon_{kl}^1=\alpha_{kl}^{-1}(\varphi)\Big((\overline {g_{kl}^1}+\hat g^1_{kl}(\varphi))\e+\hfill \\
  \hfill +(\hat f_{kl}^1(\varphi)+\Gamma_{kl}^1(t,\varphi,0)+
	\beta_{kl}(\varphi)(\overline {g_{kl}^1}+\hat g^1_{kl}(\varphi))^2)\e^2\Big) + O(\e^3).\end{matrix} $$ 

{\bf Theorem 1.} \ {\it For each $b_{kl}^*,$ for each sufficiently small  $\e>0,$
system $(\ref{sv})$ has the continuous continuously differentiable, $T$-periodic in $t,$ 
and $\omega^*$-periodic in $\varphi$ invariant surface 
\begin{equation}\label{t1} x=C(\varphi)+\Upsilon_{kl}^\nu(t,\varphi,\e)(C(\varphi)-k),\ \ 
  y=S(\varphi)+\Upsilon_{kl}^\nu(t,\varphi,\e)(S(\varphi)-l).
\end{equation}	
When $\varphi=0,$ this surface passes through a small neighborhood of the point $(b_{kl}^*,l)$  
and is homeomorphic to a two-dimensional torus, provided that the time is factored by the period. 
Here $b_{kl}^*$ is an  admissible parameter from Definition \,1, surface $\Upsilon_{kl}^\nu$ from $(\ref{surr}),$ 
$(C(\varphi),S(\varphi))$ is an  $\omega^*=\omega_{b_{kl}^*}$-periodic solution of system $(\ref{cs})$ 
with the initial values $C(0)=b_{kl}^*,\,S(0)=l.$ }

The surface \eqref{t1} is obtained by substituting the invariant surface \eqref{surr} into the change \eqref{pz}.

\section{The application of the obtained results.} 

\subsection{ The analysis of the analytic generating equation.}
Suppose that for system \eqref{sv}
\begin{equation} \label{xya} 
X(t,x,y,0)=\sum_{m,n=0}^\infty X^{(m,n)}(t)x^my^n,\ \ Y(t,x,y,0)=\sum_{m,n=0}^\infty Y^{(m,n)}x^my^n,
\end{equation}
where the power series are  absolutely convergent     uniformly in $t$   in an   open set $G^0=\{(t,x,y)\,|\, t\in\mathbb{R}^1,\ |x|<x_*,\,|y|<y_*\}$. 
Assume also that the coefficients are  real, continuous, and  $T$-periodic in $t$ coefficients.

Then, for any value of the parameter $b_{kl}$ from \eqref{nd}, that satisfies additional conditions \eqref{gamma*},\,\eqref{zig} and  \eqref{ud}, the function $R_{kl}^o$ from \eqref{rio}  takes the form
$$R_{kl}^o(t,\varphi)=\sum_{m,n=0}^\infty \left(Y^{(m,n)}(t)C'(\varphi)-X^{(m,n)}(t)S'(\varphi)\right)C^m(\varphi)S^n(\varphi).$$
%
Therefore, the left-hand side of generating equation \eqref{pu} with $\omega=\omega_{b_{kl}}$ takes the form
$$\overline {R_{kl}^o}(b_{kl})={1\over T\omega}\sum_{m=0}^\infty\sum_{n=0}^\infty
\left(\overline{Y^{(m,n)}}\int_0^{\omega} C^mS^nC'\,d\varphi-
      \overline{X^{(m,n)}}\int_0^{\omega} C^mS^nS'\,d\varphi\right).$$
When $n=0,$ the first integral is equal to zero and, when $m=0,$ the second one is equal to zero.

Integrating the  identity $(C^mS^{n+1})'=mC^{m-1}S^{n+1}C'+(n+1)C^mS^nS'$
with respect to the period and collecting coefficients, we get
$$\overline {R_{kl}^o}(b_{kl})={1\over T\omega}\sum_{m=0}^\infty\sum_{n=1}^\infty
  \!\left({m+1\over n}\overline{X^{(m+1,n-1)}}+\overline{Y^{(m,n)}}\right)I_{kl}^{m,n},\  
$$	
where $\displaystyle I_{kl}^{mn}=\int_0^{\omega} C^m(\varphi)S^n(\varphi)C'(\varphi)\,d\varphi.$
Notice, that for each of the three classes, this formula can be simplified as  follows:
\begin{multline}\label{rfinal}
 \overline {R_{00}}(b_{00})={4\over T\omega_{b_{00}}}\sum_{m,n=0}^\infty P^{(2m+1,2n+1)}J_{00}^{mn},\\ 
\overline {R_{k0}}(b_{k0})={2\over T\omega_{b_{k0}}}\sum_{m,n=0}^\infty  P^{(m,2n+1)}J_{k0}^{mn}, 
\overline {R_{kl}}(b_{kl})={1\over T\omega_{b_{kl}}}\sum_{m,n=0}^\infty P^{(m,n+1)}J_{kl}^{mn}\ \ (k,l=\pm 1),
\end{multline}
where \,$\displaystyle P^{(m,n)}={m+1\over n}\overline{X^{(m+1,n-1)}}+\overline{Y^{(m,n)}}$ \ end \ 
$\displaystyle J_{00}^{mn}=\left\{
\int_0^{ r_{0e}^1} \eta^{2m}(S^+(\eta^2))^{2n+1}\,d\eta+\right.$ \\
$\displaystyle \left.\int_{r_{0e}^1}^{b_{00}^e} \eta^{2m}(S^-(\eta^2))^{2n+1}\,d\eta \ \hbox { for }\ b_{00}^e\in (r_e,r_M),
\int_{b_{00}^i}^0 \eta^{2m}(S^-(\eta^2))^{2n+1}d\eta\ \hbox{ for }\ b_{00}^i\in (0,r_i)\right\},$ \\
$kJ_{k0}^{mn}=$ $\displaystyle \int_{k l_1^0}^{k l_1^1}\eta^m(S^-(\eta^2))^{2n+1}d\eta+ 
\int_{k l_1^1}^{k r_1^1}\eta^m(S^+(\eta^2))^{2n+1}\,d\eta+\int_{k r_1^1}^{k b_{k0}}\eta^m(S^-(\eta^2))^{2n+1}\,d\eta,$ \\
$\displaystyle kJ_{kl}^{mn}=\int_{k l_2^1}^{kb_{kl}}\eta^m((S^+(\eta^2))^{n+1}-(S^-(\eta^2))^{n+1})\,d\eta;$ \ 
$S^\pm(\eta^2)$ is from  \eqref{Spm}.

For example, let us deduce the first formula from $(\ref{rfinal})$ taking into account that the  constants are introduced in \eqref{exc}.

For the class $0_e],$ when $b_{00}^e\in (r_e,r_M),$ the  movement along a closed orbit is clockwise 
 and the function $S^2(\varphi)-1$
changes its sign in each quadrant when $|C(\varphi)|=r_{0e}^1.$ 
Therefore,
$\displaystyle I_{00}^{mn}=\int_{b_{00}^e}^{r_{0e}^1}\eta^m(-S^-(\eta^2))^n\,d\eta+\int_{ r_{0e}^1}^0 \eta^m(-S^+(\eta^2))^n\,d\eta+$ 
$\displaystyle \int_0^{-r_{0e}^1}\eta^m(-S^+(\eta^2))^n\,d\eta+\int_{-r_{0e}^1}^{-b_{00}^e}\eta^m(-S^-(\eta^2))^n\,d\eta+$ 
$\displaystyle \int_{-b_{00}^e}^{-r_{0e}^1}\eta^m( S^-(\eta^2))^n\,d\eta+\int_{-r_{0e}^1}^0 \eta^m( S^+(\eta^2))^n\,d\eta+I_e,$  
where $\displaystyle I_e=\int_0^{ r_{0e}^1}\eta^m( S^+(\eta^2))^n\,d\eta+\int_{ r_{0e}^1}^{ b_{00}^e}\eta^m( S^-(\eta^2))^n\,d\eta.$

It is obvious that $I_{00}^{mn}=4I_e,$ when $m$ is even and $n$ is odd. 
In   the other cases  $I_{00}^{mn}=0.$ 


For the  class $0_m],$ when $b_{00}^i\in (0,r_i),$ the movement along a closed orbit is counterclockwise. Therefore, 
$$
 I_{00}^{mn}=I_i+\int_0^{-b_{00}^i}\eta^m( S^-(\eta^2))^n\,d\eta+\int_{-b_{00}^i}^0\eta^m(-S^-(\eta^2))^n\,d\eta+
\int_0^{b_{00}^i} \eta^m(-S^-(\eta^2))^n\,d\eta,$$
 where $I_i=\int_{b_{00}^i}^0 \eta^m( S^-(\eta^2))^n\,d\eta.$
Then, $I_{00}^{mn}=4I_i,$ if $m$ is even and $n$ is odd. In the other cases $I_{00}^{mn}=0.$

\subsection{ The practical application of the obtained results.}

As an example, we investigate the  
$T$-periodic in $t$ system \eqref{sv} with $\gamma=0.5,$ 
$M_x,M_y\ge2.$  The functions $X(t,x,y,0)$ and $Y(t,x,y,0)$ take the  form \eqref{xya} and   
\begin{equation}\label{coeff}
\begin{array}{c} 
\forall \, q, s\in\mathbb{Z}_+\!:\ \ \overline {X^{(q,s)}}=0,\ \overline {Y^{(q,s)}}=0,\\ 
\hbox{ except }\ \overline {Y^{(0,1)}}=4.57,\ \overline {Y^{(0,3)}}=-1.66,\ \overline {Y^{(1,1)}}=-0.855,\ \overline {Y^{(0,2)}}=0.513.
\end{array}
\end{equation}

By  \eqref{exc}  for this  system  \,$r_i=2^{-1/4}(\sqrt 2-1)\approx 0.348,$ $r_e=r_\gamma=2^{-1/4}(\sqrt 2+1)\approx 1.306;$ 
according to \eqref{nd} $r_M=\sqrt 3>1.732$ (since the  constant is from Remark 1 $M=4.5$). 
Therefore, by  \eqref{nd} for the class 0]  $b_{00}^i\in (0,0.348),\ b_{00}^e\in (1.306,1.732),$ 
for class 1], $kb_{k0}\in (1.07,1.306)$.  Moreover, the  lower bound  is given taking into account conditions \eqref{gamma*},
for the class 2] $kb_{kl}\in (1,1.306)$ (here  $k,l\in\{-1,1\}$). 

In turn, formulas \eqref{rfinal} take the form: \\
$T\omega_{b_{00}} \overline {R^o_{00}}(b_{00})=4(\overline {Y^{(0,1)}} J_{00}^{00}+\overline {Y^{(0,3)}} J_{00}^{01}),$ \\
$T\omega_{b_{k0}} \overline {R^o_{k0}}(b_{k0})=2(\overline {Y^{(0,1)}} J_{k0}^{00}+\overline {Y^{(0,3)}} J_{k0}^{01}+
					  \overline {Y^{(1,1)}} J_{k0}^{10}),$ \\
$T\omega_{b_{kl}} \overline {R^o_{kl}}(b_{kl})=\overline {Y^{(0,1)}} J_{kl}^{00}+\overline {Y^{(0,3)}} J_{kl}^{02}+
					\overline {Y^{(1,1)}} J_{kl}^{10}+\overline {Y^{(0,2)}} J_{kl}^{01}.$
					
Now, for each class  we change the  initial value of $b$ in the  given limits
by the step $10^{-2},$ calculate the values of the function $T\omega_{b_{kl}} \overline {R^o_{kl}}(b_{kl})\ (k,l\in\{0,\pm 1\})$ at each step value and observe its sign.
In case we find intervals with endpoints of a different sign, it is meaningful to reiterate the process but with a smaller step value. 

For the  given system we obtain: \\	
$T\omega_{b_{00}}\overline {R^o_{00}}(b^<_{00})>10^{-5},$ $T\omega_{b_{00}}\overline {R^o_{00}}(b^>_{00})<-10^{-5}$ 
with	$b_{00}^<=1.5,\ b_{00}^>=1.501;$ \\ 
$kT\omega_{b_{k0}}\overline {R^o_{k0}}(b^<_{k0})<-10^{-5},$ $k\omega_{b_{k0}}T\overline {R^o_{k0}}(b^>_{k0})>10^{-5}$
for	$b_{k0}^<=1.2,\ b_{k0}^>=1.202;$ \\
$klT\omega_{b_{kl}}\overline {R^o_{kl}}(b^<_{kl})<-10^{-5},$ $kl\omega_{b_{kl}}T\overline {R^o_{kl}}(b^>_{kl})>10^{-5}$
for  	$b_{kl}^<=1.2,\ b_{kl}^>=1.202.$ 
	 
Therefore, we find seven values of the parameter $b:$ 
$b_{00}^{e0}\in (1.5,1.501),$ $kb_{k0}^0,kb_{k1}^0,kb_{k,-1}^0\in (1.2,1.202)$ $(k=\pm 1).$ 
These values are the zeros of generating equation \eqref{pu} and belong to the corresponding intervals from \eqref{nd} and \eqref{gamma*}.
By Definition 1 the mentioned values are admissible if the  Siegel and dissipativity conditions hold.

In such case by Theorem 1  under conditions \eqref{coeff}, for each sufficiently small $\varepsilon>0$ 
system \eqref{sv} has seven two-dimensional invariant surfaces \eqref{t1}. 
When  $t=0,$ four of them  generate closed orbits of class 2], which enclose one of  singular points $(1,\pm 1),(-1,\pm 1)$ 
of the unperturbed system \eqref{snv}, another two surfaces are  closed orbits of class 1], which enclose one of  separatrix "eight"\  
and three singular points: $(\pm 1,1),(\pm 1,0),(\pm 1,-1).$ Finally, the seventh surface is  of the class $0^e]$ and  encloses all nine singular pints.

We have not found  the eighth invariant surface   generated by a cycle  of the class $0^i]$ which encloses only the point $(0,0),$ since 
for each $b_{00}^i\in (0,0.348)$  we have $\overline {R^o_{00}}(b_{00}^i)<0.$  

Calculations, that prove that all seven constants $L$ from dissipativity condition \eqref{ud} are not equal to zero, are not provided due to their cumbersome nature.
To perform these calculations for each initial value $b_{00},$ we interpolated solutions of the 
initial value problems $C(\varphi),S(\varphi)$ for  system \eqref{cs} with the fifth degree polynomials depending on $\varphi.$ 

Now, the only thing left to do is to notice that the Siegel condition, required for the case $\nu=0,$ 
holds for almost all vectors with integer components 
with respect to the Lebesgue measure.

\section*{Acknowledgments}

Valery Romanovski  acknowledges the financial support from the Slovenian Research Agency (research core funding  P\,1-0306 and project N\,1-0063).

\end{document}